\documentclass[draft]{amsart}

\usepackage{ifthen}

\newtheorem{anyprop}{Anyprop}[section]

\newtheorem{theorem}[anyprop]{Theorem}

\newtheorem{proposition}[anyprop]{Proposition}
\newtheorem{corollary}[anyprop]{Corollary}

\theoremstyle{definition}

\newtheorem{example}[anyprop]{Example}
\newtheorem{examples}[anyprop]{Examples}

\newtheorem{remark}[anyprop]{Remark}

\newcommand{\ZZ}{\mathbb{Z}}

\newcommand{\PP}{\mathbb{P}}

\newcommand  {\shE}     {\mathcal{E}}
\newcommand  {\shF}     {\mathcal{F}}

\newcommand  {\shS}     {\mathcal{S}}

\renewcommand  {\ker }  {\operatorname{ker}}

\renewcommand{\O}       {\mathcal{O}}

\newcommand  {\Proj}    {\operatorname{Proj}}

\newcommand  {\rk}    {\operatorname{rk}}

\newcommand  {\Syz}     {\operatorname{Syz}}

\newcommand{\komdots}{ , \ldots , }

\newcommand{\lto}{\longrightarrow}

\theoremstyle{remark}

\numberwithin{equation}{section}

\usepackage{amscd}
\usepackage{amssymb}

\begin{document}

\newboolean{brennersupplement}

\setboolean{brennersupplement}{false}

\title[Syzygy Bundles on $\PP^2$ and the Weak Lefschetz Property]
{Syzygy Bundles on $\PP^2$ and the Weak Lefschetz Property}

\author[Holger Brenner and Almar Kaid]{Holger Brenner and Almar Kaid}
\address{Department of Pure Mathematics, University of Sheffield,
  Hicks Building, Hounsfield Road, Sheffield S3 7RH, United Kingdom}
\email{H.Brenner@sheffield.ac.uk and A.Kaid@sheffield.ac.uk}


\subjclass{}



\begin{abstract}
Let $K$ be an algebraically closed field of characteristic zero and
let $I=(f_1 \komdots f_n)$ be a homogeneous $R_+$-primary ideal in
$R:=K[X,Y,Z]$. If the corresponding syzygy bundle $\Syz(f_1 \komdots
f_n)$ on the projective plane is semistable, we show that the
Artinian algebra $R/I$ has the Weak Lefschetz property if and only
if the syzygy bundle has a special generic splitting type. As a
corollary we get the result of Harima et alt., that every Artinian
complete intersection ($n=3$) has the Weak Lefschetz property.
Furthermore, we show that an almost complete intersection ($n=4$)
does not necessarily have the Weak Lefschetz property, answering
negatively a question of Migliore and Mir\'{o}-Roig. We prove that
an almost complete intersection has the Weak Lefschetz property if
the corresponding syzygy bundle is not semistable.
\end{abstract}

\maketitle

Mathematical Subject Classification (2000): primary: 13D02, 14J60,
secondary: 13C13, 13C40, 14F05.

Keywords: syzygy, semistable bundle, Grauert-M\"ulich Theorem, Weak
Lefschetz property, Artinian algebra, complete intersection, almost
complete intersection.

\section{Introduction}

Throughout this paper we denote by $R:=K[X_0\komdots X_N]$ a
polynomial ring in $N+1$ variables over an algebraically closed
field $K$ of characteristic zero. A family of $R_+$-primary
homogeneous polynomials $f_1 \komdots f_n$ of degree $d_i$ in $R$
defines on $\PP^N = \PP^N_K = \Proj R$ the short exact sequence
$$0 \lto \Syz(f_1 \komdots f_n)(m) \lto \bigoplus_{i=1}^n \O_{\PP^N}
(m-d_i) \stackrel{f_1 \komdots f_n}{\lto} \O_{\PP^N}(m) \lto 0$$
of locally free sheaves. We call the vector bundle of rank $n-1$
on the left the \emph{syzygy bundle} of the elements $f_1 \komdots
f_n$. On the other hand the elements $f_1 \komdots f_n$ define an
Artinian graded $R$-algebra $A:=R/(f_1 \komdots f_n)$, i.e. $A$ is
of the form
$$A = K \oplus A_1 \oplus \ldots \oplus A_s$$
for an integer $s \geq 0$. The algebra $A$ has the so called
\emph{Weak Lefschetz property {\rm(}WLP\,{\rm)}} if for every
general linear form $\ell \in R_1$ the multiplication maps
$$A_m \stackrel{\cdot \ell} \longrightarrow A_{m+1}$$
have maximal rank for $m=0 \komdots s-1$, i.e. these $K$-linear maps
are either injective or surjective.

In case of generic elements $f_1 \komdots f_n \in R$, the Weak
Lefschetz property is related to the study of the Fr\"oberg
conjecture. This conjecture is equivalent to the \emph{Maximal Rank
property}, i.e. the property that the multiplication maps have
maximal rank for every $d \geq 1$ and every general form $F \in
R_d$. The Fr\"oberg conjecture is only known for the cases $N=1,2$
(cf. \cite{froeberg} and \cite{anick}).

Our objective is to study for $N= 2$ the connection between the WLP
for the algebra $A$ and the semistability of the syzygy bundle
$\Syz(f_1 \komdots f_n)$. We recall that a torsion-free coherent
sheaf $\shE$ on $\PP^N$ is \emph{semistable} (in the sense of
Mumford and Takemoto) if for every coherent subsheaf $\shF \subset
\shE$ the inequality $\frac{\deg(\shF)}{\rk(\shF)} \leq
\frac{\deg(\shE)}{\rk(\shE)}$ holds (and \emph{stable} if $<$
holds), where the \emph{degree} $\deg(\shF)$ is defined as the twist
$k$ such that $(\bigwedge^{\rk(\shF)} \shF)^{**} \cong
\O_{\PP^N}(k)$ (cf. \cite{huybrechtslehn} and
\cite{okonekschneiderspindler}). By the Theorem of
\emph{Grauert-M\"ulich} (cf. \cite[Theorem 3.0.1]{huybrechtslehn} or
\cite[Corollary 1 of Theorem 2.1.4]{okonekschneiderspindler}) a
semistable vector bundle $\shE$ of rank $r$ on $\PP^N$ splits on a
generic line $L \subset \PP^N$ as $\shE|_L \cong \O_L (a_1)
\oplus \cdots \oplus \O_L (a_r)$ with $a_1 \geq \ldots \geq
a_r$ and $0 \leq a_i-a_{i+1} \leq 1$ for $i=1 \komdots r-1$. We
prove in Theorem \ref{genericsplittypeforsyz} that in case of a
semistable syzygy bundle on $\PP^2$ the algebra $A$ has $WLP$ if and
only if in the generic splitting type of $\Syz(f_1\komdots f_n)$ at
most two different twists occur. The importance of the Theorem of
Grauert-M\"ulich was already mentioned in \cite{harimamigliore}, but
only for complete intersections in $\PP^2$.

In the case of an almost complete intersection (i.e. four ideal
generators) we show in Example \ref{counterexample} that for a
semistable syzygy bundle all numerically possible splitting types do
actually exist. It follows that there are examples of almost
complete intersections where the WLP does not hold, which gives a
negative answer to a question of Migliore and Mir\'{o}-Roig (cf.
\cite[Paragraph after Question 4.2]{miglioremiroroiglefschetz}).
Furthermore, we prove in Theorem \ref{almostcompleteintcase} that in
the non semistable case an almost complete intersection has always
the Weak Lefschetz property.

\section{Generic splitting type of syzygy bundles and the Weak
Lefschetz property}

We start with the following cohomological observation.

\begin{proposition}\label{h1gleichstufe}
Let $R=K[X_0,\ldots,X_N]$, $N \geq 2$, and let
$I=(f_1,\ldots,f_n)\subseteq R$ be an $R_+$-primary homogeneous
ideal. Then we have
$$A_m=H^1(\PP^N,\Syz(f_1\komdots f_n)(m))$$ for every graded
component $A_m$ of $A:=R/I$, $m \in \ZZ$.
\end{proposition}

\begin{proof}
Since $H^1(\PP^N,\O_{\PP^N}(m))=0$ for all $m$ and $N \geq 2$, we
derive from the presenting sequence of $\shS(m):=\Syz(f_1\komdots
f_n)(m)$ the exact cohomology sequence
$$\bigoplus_{i=1}^n H^0(\PP^N,\O_{\PP^N}(m-d_i)) \stackrel{f_1 \komdots f_n} \lto H^0(\PP^N,\O_{\PP^N}(m)) \lto H^1(\PP^N,\shS(m))
\lto 0.$$ Now the claim follows immediately, since
$H^0(\PP^N,\O_{\PP^N}(m)) = R_m$.
\end{proof}

We restrict now to three variables and write $R=K[X,Y,Z]$. As usual
let $f_1 \komdots f_n$ be homogeneous $R_+$-primary elements and let
$0 \neq \ell \in R_1$ be a linear form. Thus $\ell$ defines an exact
sequence
$$0 \lto \O_{\PP^2} \stackrel{\cdot \ell}\lto \O_{\PP^2}(1) \lto
\O_L(1) \lto 0,$$ where the map on the right is the restriction to
the line $L = V_+(\ell) \subset \PP^2$ defined by $\ell$. Since
the syzygy bundle $\shS(m):=\Syz(f_1 \komdots f_n)(m)$ is locally
free this yields the short exact sequence
$$0 \lto \shS(m) \stackrel{\cdot \ell}\lto \shS(m+1)
\lto \shS_L(m+1) \lto 0$$ (where $\shS_L:=\shS|_L$), and from this
we derive the long exact sequence in cohomology
\begin{eqnarray*}
0 &\lto& H^0(\PP^2,\shS(m)) \stackrel{\cdot \ell}\lto
H^0(\PP^2,\shS(m+1)) \lto H^0(L,\shS_L(m+1))\\
&\stackrel{\delta}\lto& A_m = H^1(\PP^2,\shS(m)) \stackrel{\cdot
\ell} \lto A_{m+1}= H^1(\PP^2,\shS(m+1))\\ &\lto&
H^1(L,\shS_L(m+1))\stackrel{\delta}\lto H^2(\PP^2,\shS(m)).
\end{eqnarray*}

Hence the Artinian algebra $A = R/(f_1 \komdots f_n)$ has the Weak
Lefschetz property if and only if for every generic line $L
\subset \PP^2$ the map $H^1(\PP^2,\shS(m)) \rightarrow
H^1(\PP^2,\shS(m+1))$ is either injective or surjective. The
injectivity is equivalent to the surjectivity of the restriction map
$H^0(\PP^2,\shS(m+1)) \rightarrow H^0(L,\shS_L(m+1))$ and holds
in particular when $H^0(L,\shS_L(m+1)) = 0$. The surjectivity is
equivalent to the injectivity of $H^1(L,\shS_L(m+1))
\stackrel{\delta}\rightarrow H^2(\PP^2,\shS(m))$ and holds in
particular when $H^1(L,\shS_L(m+1))=0$. From now on we denote the
map $A_m \stackrel{\cdot \ell} \rightarrow A_{m+1}$ by $\mu_\ell$.

The following theorem relates the generic splitting type of a
semistable syzygy bundle with the Weak Lefschetz property of the
Artinian algebra $A$.

\begin{theorem}\label{genericsplittypeforsyz}
Let $f_1 \komdots f_n$ be $R_+$-primary homogeneous polynomials in
$R=K[X,Y,Z]$ such that their syzygy bundle $\shS:=\Syz(f_1 \komdots
f_n)$ is semi\-stable on $\PP^2$. Then the following holds:
\begin{enumerate}
\item If the restriction of $\shS$ splits on a generic line $L$ as
$$\shS|_L \cong \bigoplus_{i=1}^{s}\O_L(a+1) \oplus
\bigoplus_{i=s+1}^{n-1} \O_L(a),$$ then $A=R/(f_1 \komdots f_n)$
has WLP. \item If the restriction of $\shS$ splits on a generic
line $L$ as
$$\shS|_L \cong \O_L(a_1) \oplus \ldots \oplus \O_L(a_{n-1})$$
with $a_1 \geq a_2 \geq \ldots \geq a_{n-1}$ and $a_1 - a_{n-1}
\geq 2$, then $A = R/(f_1\komdots f_n)$ has not WLP.
\end{enumerate}
\end{theorem}

\begin{proof}
Let $\ell$ be the general linear form defining a generic line $L
\subset \PP^2$. To prove the first part, according to Proposition
\ref{h1gleichstufe}, we have to show that the multiplication map
$\mu_\ell$ is either surjective or injective for every $m \in \ZZ$.
So we consider the long exact sequence in cohomology mentioned
above. Firstly, we assume $m < -a-2$. But then $H^0(L,\shS_L(m+1))=
\bigoplus_{i=1}^s H^0(L,\O_L(m+a+1)) \oplus \bigoplus_{i=s+1}^{n-1}
H^0(L,\O_L(m+a))=0$ and hence $\mu_\ell$ is injective. Now let $m
\geq -a-2$. Then Serre duality yields $H^1(L,\shS_L(m+1)) \cong
H^0(L,\shS_L^*(-m-3))^*$. Since the dual bundle $\shS^*(-m-3)$
splits on $L$ as $\bigoplus_{i=1}^{s}\O_L(-a-m-4) \oplus
\bigoplus_{i=s+1}^{n-1} \O_L(-a-m-3)$, it has no non-trivial
sections on $L$. Hence the map $\mu_\ell$ is onto.

For the proof of part two we observe that for the degrees
$\deg(\shS(-a_1)) < 0$ and $\deg(\shS(-a_{n-1})) > 0$ holds. Since
$a_1 - a_{n-1} \geq 2$, we can find an $m \in \ZZ$ with $-a_{n-1}-2
\geq m+1 \geq -a_1$ and such that $\deg(\shS(m+1)) < 0$ and
$\deg(\shS(m+3))
> 0$. For this $m$ we have $H^0(\PP^2,\shS(m+1))=0$, since the syzygy bundle
$\shS(m+1)$ is semistable on $\PP^2$. Because of $a_1 + m + 1 \geq
0$ we also have
$H^0(L,\shS_L(m+1))=\bigoplus_{i=1}^{n-1}H^0(L,\O_L(a_i+m+1)) \neq
0$. Hence the map $H^0(\PP^2,\shS(m+1)) \rightarrow
H^0(L,\shS_L(m+1))$ is not surjective and thus $\mu_\ell$ is not
injective in degree $m$. Now we show that $\mu_\ell$ is not
surjective either in this particular degree $m$. Serre duality
yields $H^1(L,\shS_L(m+1)) \cong H^0(L,\shS_L^*(-m-3))^*$ and
$H^2(\PP^2,\shS(m)) \cong H^0(\PP^2,\shS^*(-m-3))^*$. Since
$\shS^*(-m-3)$ is semistable and $\deg(\shS^*(-m-3))<0$ we have
$H^0(\PP^2,\shS^*(-m-3))=0$. Now we conclude as above that
$H^0(L,\shS^*_L(-m-3))=\bigoplus_{i=1}^{n-1}H^0(L,\O_L(-a_i-m-3))
\neq 0$ because $-a_{n-1}-m-3 \geq 0$ is equivalent to $m+1 \leq
-a_{n-1}-2$. So $\mu_\ell$ is not surjective.
\end{proof}

\begin{remark}
Since we have not used any particular properties of syzygy
bundles, Theorem \ref{genericsplittypeforsyz} can be generalized
for arbitrary vector bundles on $\PP^2$ if one translates the Weak
Lefschetz property for a vector bundle $\shE$ into the property
that the multiplication map $H^1(\PP^2,\shE(m)) \rightarrow
H^1(\PP^2,\shE(m+1))$ induced by a general linear form has maximal
rank.
\end{remark}

We can now prove \cite[Theorem 2.3]{harimamigliore} easily.

\begin{corollary}\label{completeintcase}
Every Artinian complete intersection in $K[X,Y,Z]$ has the Weak
Lefschetz property.
\end{corollary}

\begin{proof}
Let $f_1,f_2,f_3 \in R$ be an Artinian complete intersection and let
$\ell \in R_1$  be a generic linear form. Firstly, we consider the
case that the corresponding syzygy bundle $\shS:=\Syz(f_1,f_2,f_3)$
is semistable. Since $\shS$ is a $2$-bundle, its restriction to the
generic line $L$ defined by $\ell$ splits by the Theorem of
Grauert-M\"ulich as $\O_L(a_1) \oplus \O_L(a_2)$ with $a_1 \geq a_2$
and $0 \leq a_1-a_2 \leq 1$. Hence by Theorem
\ref{genericsplittypeforsyz}(1) the algebra $R/(f_1,f_2,f_3)$ has
WLP.

Now suppose $\shS$ is not semistable. Since we can pass over to the
reflexive hull, the Harder Narasimhan filtration (cf.
\cite[Definition 1.3.2]{huybrechtslehn}) of $\shS$ looks like $0
\subset \O_{\PP^2}(a) \subset \shS$ with $a \in \ZZ$ (cf.
\cite[Lemma 1.1.10]{okonekschneiderspindler}). The quotient
$\shF:=\shS/\O_{\PP^2}(a)$ is a torsion-free sheaf, which is outside
codimension $2$ isomorphic to its bidual $\shF^{**}$. This bidual is
reflexive (cf. \cite[Lemma 24.2]{schejastorchverzweigung}), hence
locally free on $\PP^2$ (cf. \cite[Lemma
1.1.10]{okonekschneiderspindler}), i.e. $\shF^{**} \cong
\O_{\PP^2}(b)$, for $b \in \ZZ$ and $a> b$. Hence we have $\shS|_L
\cong \O_L(a) \oplus \O_L(b)$ for a generic line $L$. If $m < -a-1$
then $H^0(L,\shS_L(m+1)) = 0$ and the multiplication map $\mu_\ell$
is injective. Now, assume that $m \geq -a-1$. We apply Serre duality
and get
\begin{eqnarray*}
H^1(L,\shS_L(m+1))\! &=&\! H^0(L,\shS_L^*(-m-3))^*\\ \! &=&\!
H^0(L,\O_L(-a-m-3))^* \oplus H^0(L,\O_L(-b-m-3))^*\\
\!&=&\! H^0(L,\O_L(-b-m-3))^* = H^1(L,\O_L(b+m+1)).
\end{eqnarray*}
The sheaf morphism $\shS \rightarrow \shF \rightarrow \shF^{**}
\cong \O_{\PP^2}(b)$ induces a map $H^2(\PP^2,\shS) \rightarrow
H^2(\PP^2,\O_{\PP^2}(b))$. Hence we have by the functoriality of the
connecting homomorphism (cf. \cite[Theorem III.1.1.A(d)]{haralg}) a
commutative diagram
$$\begin{CD}
H^1(L,\shS_L(m+1)) @>\delta>> H^2(\PP^2,\shS(m))\\
@V\cong VV                       @VVV\\
H^1(L,\O_L(b+m+1)) @>\delta>> H^2(\PP^2,\O_{\PP^2}(b+m)),
\end{CD}$$
where the bottom map is injective by its explicit description.
Hence the map $H^1(L,\shS_L(m+1)) \stackrel{\delta} \rightarrow
H^2(\PP^2,\shS(m))$ is injective as well.
\end{proof}

\begin{corollary}
Let $f_1 \komdots f_n$, $n \geq 3$, be generic forms in
$R=K[X,Y,Z]$ such that their syzygy bundle is semistable. Then
$\Syz(f_1 \komdots f_n)$ splits on a generic line $L$ as
$$\Syz(f_1 \komdots f_n)|_L \cong \bigoplus_{i=1}^{s}\O_L(a+1)
\oplus \bigoplus_{i=s+1}^{n-1}\O_L(a).$$
\ifthenelse{\boolean{brennersupplement}}{In particular, the syzygy
bundle for generic forms of constant degree splits like this.}{In
particular, if $f_1 \komdots f_n$ are generic polynomials of
constant degree $d$ and $3\leq n \leq 3d$ then the corresponding
syzygy bundle has this splitting type.}
\end{corollary}

\begin{proof}
Since Anick proved in \cite[Corollary 4.14]{anick} that every
ideal of generic forms $f_1 \komdots f_n$, $n \geq 3$, in
$K[X,Y,Z]$ has the Weak Lefschetz property (in fact Anick proved a
much stronger result), this follows immediately from Theorem
\ref{genericsplittypeforsyz}(2).\ifthenelse{\boolean{brennersupplement}}{The
supplement follows from \cite{brennermonomial}, where it was
proved that the syzygy bundle for generic forms of a fixed degree
is semistable on $\PP^2$.}{ The supplement follows from
\cite[Theorem A.1]{heinsyzygystable}.}
\end{proof}

\section{Almost complete intersections}

In \cite[Paragraph after Question 4.2]{miglioremiroroiglefschetz},
Migliore and Mir\'{o}-Roig ask whether \emph{every} almost complete
intersection (i.e. four ideal generators) in $K[X,Y,Z]$ has the Weak
Lefschetz property. The following easy example gives via Theorem
\ref{genericsplittypeforsyz} a negative answer to this question.

\begin{example}\label{counterexample}
We consider the monomial almost complete intersection generated by
$X^3,Y^3,Z^3,XYZ$ in $K[X,Y,Z]$. The corresponding syzygy bundle
$\Syz(X^3,Y^3,Z^3,XYZ)$ is semistable by \cite[Corollary
3.6]{brennerlookingstable}. We compute its restriction to a line
$L$ given by $Z = uX + vY$ with arbitrary coefficients $u,v \in
K$. We have
$$Z^3|_L = u^3X^3+v^3Y^3 + 3uv(uX^2Y+vXY^2)~\mbox{  and  }~(XYZ)|_L=uX^2Y+vXY^2.$$
For $u,v \neq 0$, and in particular for generic $u,v$, this gives
immediately the non-trivial syzygy
$$u^3X^3+v^3Y^3+3uv(uX^2Y+vXY^2)-u^3X^3-v^3Y^3-3uv(uX^2Y+vXY^2),$$
which yields a non-trivial global section in
$\Syz(X^3,Y^3,Z^3,XYZ)(3)|_L$. Since $\Syz(X^3,Y^3,Z^3,XYZ)(3)$
has degree $-3<0$, this section does not come from $\PP^2$.
Moreover,
$$\Syz(X^3,Y^3,Z^3,XYZ)(4)|_L \cong \O_L(1) \oplus
\O_L \oplus \O_L(-1),$$ hence by Theorem
\ref{genericsplittypeforsyz} the Artinian algebra $$A =
K[X,Y,Z]/(X^3,Y^3,Z^3,XYZ)$$ has not the Weak Lefschetz property.
Indeed, the map $A_2 \rightarrow A_3$ given by a generic linear form
is neither injective nor surjective since $v^2X^2 +
u^2Y^2+Z^2-uvXY-vXZ-uYZ$ is in the kernel for every generic linear
form $\ell = uX + vY + Z$ and $\dim_K A_2 = 6 = \dim_K A_3$.

The following reasoning shows that this is the only counterexample
in degree $3$ containing the monomials $X^3,Y^3,Z^3$. So we consider
the monomials $X^3,Y^3,Z^3$ and a forth homogeneous polynomial $f =
\sum_{|\nu|=3} a_\nu X^\nu$ of degree $3$. If we use again $Z=uX+vY$
to restrict the corresponding syzygy bundle to a generic line, we
have to compute the coefficients $c_{(1,2)}$ and $c_{(2,1)}$ of the
monomials $XY^2$ and $X^2Y$ in $f$ restricted to $L$. These are
$$c_{(1,2)} = a_{(1,2,0)} + a_{(0,2,1)}u + 2a_{(0,1,2)}uv +
a_{(1,1,1)}v + a_{(1,0,2)}v^2 + 3 a_{(0,0,3)} u v^2$$ and
$$c_{(2,1)} = a_{(2,1,0)} + a_{(0,1,2)}u^2 + a_{(2,0,1)}v +
a_{(1,1,1)}u + 2a_{(1,0,2)}uv + 3 a_{(0,0,3)}u^2v.$$ The algebra $A
= K[X,Y,Z]/(X^3,Y^3,Z^3,f)$ has not the Weak Lefschetz property if
and only if there exists a non-trivial global section of $S|_L(3)$,
and this is true if $c_{(2,1)}X^2Y + c_{(1,2)}XY^2$ is a multiple of
$uX^2Y+vXY^2$, i.e. $c_{(2,1)}X^2Y + c_{(1,2)}XY^2 = t (
uX^2Y+vXY^2)$ for some $t \in K$. This means $c_{(2,1)} = t u$ and
$c_{(1,2)}=tv$ and gives the condition $vc_{(2,1)} - uc_{(1,2)} =0$.
We have $vc_{(2,1)} - uc_{(1,2)} =$
$$va_{(2,1,0)}-ua_{(1,2,0)}+v^2a_{(2,0,1)}-u^2a_{(0,2,1)}-u^2va_{(0,1,2)}+uv^2a_{(1,0,2)}.$$
If we consider the right hand side of this equation as a polynomial
in $K[u,v]$ and assume that at least one of the coefficients is not
zero, then there exists also values in $K$ where this polynomial
does not vanish ($K$ is an infinite field). Hence we see that if $f
\notin (X^3,Y^3,Z^3,XYZ)$ this condition can not hold for all $u,v$
and therefore the algebra $A$ has the Weak Lefschetz property.
\end{example}

\begin{remark}
In \cite[Question 4.2]{miglioremiroroiglefschetz}, Migliore and
Mir\'{o}-Roig asked: ``For any integer $n \geq 3$, find the maximum
number $A(n)$ (if it exists) such that \emph{every} Artinian ideal
$I \subset k[x_1 \komdots x_n]$ with $\mu(I) \leq A(n)$ has the Weak
Lefschetz property (where $\mu(I)$ is the minimum number of
generators of $I$).'' They show in \cite[Example
4.2]{miglioremiroroiglefschetz} that the Artinian ideal
$I:=(X^2,Y^2,Z^2,XY,XZ)$ has not the Weak Lefschetz property, so
$A(3) \leq 4$. Since by Corollary \ref{completeintcase} every
complete intersection in $K[X,Y,Z]$ has the Weak Lefschetz property,
Example \ref{counterexample} proves that $A(3)=3$. Furthermore,
Example \ref{counterexample} shows also that
$$XYZ \in \ker[(R/(X^3,Y^3,Z^3))_3 \lto
(R/(\ell,X^3,Y^3,Z^3))_3]$$ for all linear forms $\ell \in R_1$.
Hence, the last hypothesis in \cite[Proposition
5.5]{miglioremiroroiglefschetz} does not always hold.
\end{remark}

For an almost complete intersection in $K[X,Y,Z]$ we can prove the
following theorem.

\begin{theorem}\label{almostcompleteintcase}
Let $f_1,f_2,f_3,f_4$ be $R_+$-primary homogeneous polynomials in
$R=K[X,Y,Z]$ such that their syzygy bundle is not semistable. Then
the algebra $R/(f_1,f_2,f_3,f_4)$ has the Weak Lefschetz property.
\end{theorem}

\begin{proof}
We consider the Harder-Narasimhan filtration of the vector bundle
$\shS:=\Syz(f_1,f_2,f_3,f_4)$. As in the proof of Corollary
\ref{completeintcase} we can pass over to the reflexive hull and
since reflexive sheaves on $\PP^2$ are locally free (cf. \cite[Lemma
1.1.10]{okonekschneiderspindler}), we can assume that all subsheaves
in the HN-filtration are vector bundles. Further we fix a generic
line $L=V_+(\ell)$, where $\ell \in R_1$ is a general linear form.

Firstly, we treat the case that the destabilizing subbundle is a
line bundle $\O_{\PP^2}(a_1)$ for some $a_1 \in \ZZ$, i.e. we have
an exact sequence $0 \rightarrow \O_{\PP^2}(a_1)\rightarrow \shS
\rightarrow \shE \rightarrow 0$ with a semistable torsion-free sheaf
$\shE$ of rank $2$. We apply the Theorem of Grauert-M\"ulich to
$\shE^{**}$, which is locally free and outside codimension $2$
isomorphic to $\shE$. Therefore we get $\shE|_L \cong \O_L(a_2)
\oplus \O_L(a_3)$ with $a_2 \geq a_3$ and $0 \leq a_2-a_3 \leq 1$.
Since the slopes in the HN-filtration decrease strictly we have
$2a_1>a_2 + a_3 \geq a_2 + a_2 -1 = 2a_2-1$ and therefore $a_1 \geq
a_2 \geq a_3$. In particular $\shS|_L \cong \O_L(a_1) \oplus
\O_L(a_2) \oplus \O_L(a_3)$. For $m <-a_2-1$ we have
$H^0(L,\O_L(a_2+m+1)\oplus \O_L(a_3+m+1))=0$ (for $m<-a_1-1$ we have
even $H^0(L,\shS_L(m+1))=0)$. Therefore the injectivity of
$\mu_\ell$ follows from the diagram
$$\begin{CD}
H^0(\PP^2,\O_{\PP^2}(a_1+m+1)) @>>> H^0(L,\O_L(a_1+m+1))\\
@VVV                                  @VV \cong V\\
H^0(\PP^2,\shS(m+1)) @>>> H^0(L,\shS_L(m+1)).
\end{CD}$$
For $m \geq -a_2-1$ we get $H^1(L,\shS_L(m+1))=
H^0(L,\shS_L^*(-m-3))^*=0$, since $-a_3-m-3\leq
1-a_2-m-3=-a_2-m-2<0$ holds. Hence $\mu_\ell$ is surjective.

Now suppose that $\shS$ is destabilized by a semistable vector
bundle $\shE$ of rank $2$. (This case corresponds essentially to
\cite[Proposition 5.2]{miglioremiroroiglefschetz}, where one of the
ideal generators has a sufficiently large degree.) Here we apply the
Theorem of Grauert-M\"ulich to $\shE$ and get $\shE|_L \cong
\O_L(a_1)\oplus \O_L(a_2)$ with $a_1 \geq a_2$ and $0\leq
a_1-a_2\leq 1$. The quotient $\shF:=\shS/\shE$ is outside
codimension $2$ isomorphic to $\O_{\PP^2}(a_3)$ with $a_3 \in \ZZ$
and since $\shF$ is a quotient in the HN-filtration we get again
$a_1 \geq a_2 \geq a_3$ and $S|_L \cong \O_L(a_1) \oplus \O_L(a_2)
\oplus \O_L(a_3)$. Hence for $m < -a_1-1$ we have
$H^0(L,\shS_L(m+1))=0$ and thus $\mu_\ell$ is injective. We treat
the case $m \geq -a_1-1$ similar to the analog situation in the
proof of Corollary \ref{completeintcase}. By Serre duality we get
$H^1(L,\shS_L(m+1)) = H^0(L,\shS^*_L(-m-3))^*=
H^0(L,\O_L(-a_3-m-3))^*=H^1(L,\O_L(a_3+m+1))$, since
$H^0(L,\O_L(-a_1-m-3) \oplus \O_L(-a_2-m-3))=0$. The map $\shS
\rightarrow \shF \rightarrow \shF^{**} \cong \O_{\PP^2}(a_3)$
induces a map $H^2(\PP^2,\shS) \rightarrow
H^2(\PP^2,\O_{\PP^2}(a_3))$ between the second cohomology groups.
Therefore the injectivity of the map $H^1(L,\shS_L(m+1))
\stackrel{\delta}\rightarrow H^2(\PP^2,\shS(m))$ follows from the
injectivity of the map $H^1(L,\O_L(a_3+m+1)) \stackrel{\delta}
\rightarrow H^2(\PP^2,\O_{\PP^2}(a_3+m))$. Hence the map $\mu_\ell$
is surjective and the Artinian algebra $R/(f_1,f_2,f_3,f_4)$ has
WLP.

Finally we consider a HN-filtration of the form $0 \subset
\O_{\PP^2}(a_1) \subset \shE \subset \shS$ with a vector bundle
$\shE$ of rank $2$. Since the quotients of this filtration are
torsion-free, we have outside codimension $2$ the identifications
$\shE/\O_{\PP^2}(a_1) \cong \O_{\PP^2}(a_2)$ and $\shS/\shE \cong
\O_{\PP^2}(a_3)$ with $a_1 > a_2> a_3$. Therefore, we have on the
generic line $L$ the splitting $\shS|_L \cong \O_{\PP^2}(a_1) \oplus
\O_{\PP^2}(a_2) \oplus \O_{\PP^2}(a_3)$. For $m < -a_2-1$ we have
$H^0(L,\O_L(a_2+m+1) \oplus \O_L(a_3+m+1)) = 0$, i.e.
$H^0(L,\shS_L(m+1)) = H^0(L,\O_L(a_1))$. Since $0 \rightarrow
\O(a_1) \rightarrow \shS$ the map $H^0(\PP^2,\shS(m+1)) \rightarrow
H^0(L,\shS_L(m+1))$ is onto. Now, let $m \geq -a_2-1$. Here we have
$H^1(L,\shS_L(m+1)) = H^1(L,\O_L(a_3+m+1))$ and $\shS \rightarrow
(\shS/\shE)^{**} \cong \O_{\PP^2}(a_3)$. Hence the injectivity of
the map $H^1(L,\shS_L(m+1)) \stackrel{\delta}\rightarrow
H^2(\PP^2,\shS(m))$ follows again from the injectivity of
$H^1(L,\O_L(a_3+m+1)) \stackrel{\delta} \rightarrow
H^2(\PP^2,\O_{\PP^2}(a_3+m))$ and we are done.
\end{proof}

We want to give examples which show that all the three possible
Harder-Narasimhan filtrations mentioned in the proof of Theorem
\ref{almostcompleteintcase} can appear. It is even possible to
provide monomial examples for all these cases. For degree
computations of syzygy bundles see \cite[Lemma
2.1]{brennerlookingstable}.

\begin{examples}
The syzygy bundle for the monomials $X^4,Y^4,Z^4,X^3Y$ has degree
$-16$, hence its slope equals $-16/3\approx-5.33$. This vector
bundle is not semistable because the subsheaf $\Syz(X^4,X^3Y)$ has
degree $-5>-5.33$. Since the degrees are constant, there are no
maps into line bundles which contradict the semistability.
Therefore
$$\O_{\PP^2}(-5)\cong \Syz(X^4,X^3Y) \subset
\Syz(X^4,Y^4,Z^4,X^3Y)$$ constitutes the HN-filtration of the syzygy
bundle. The generic splitting type is $\O_L(-5) \oplus \O_L(-5)
\oplus \O_L(-6)$.

To give an example for the second type of HN-filtration consider the
monomials $X^4,Y^4,Z^4,X^3Y^3Z^3$. The corresponding syzygy bundle
has the slope $-21/3=-7$ and it is destabilized by the semistable
subsheaf $\Syz(X^4,Y^4,Z^4)$ of slope $-12/2=-6$. Hence we have
found the HN-filtration of the bundle $\Syz(X^4,Y^4,Z^4,X^3Y^3Z^3)$.
The generic splitting type is $\O_L(-6) \oplus \O_L(-6) \oplus(-9)$.

For the third type of HN-filtration we look at the family
$X^2,Y^4,Z^7,XY$. The HN-filtration of $\Syz(X^2,Y^4,Z^7,XY)$ is
$$0 \subset \O_{\PP^2}(-3) \cong \Syz(X^2,XY) \subset
\Syz(X^2,XY,Y^4) \subset \Syz(X^2,Y^4,Z^7,XY),$$ since the
quotients have rank one and degrees $-3 > -5 > -7$. Accordingly,
the generic splitting type is $\O_L(-3) \oplus \O_L(-5) \oplus
\O_L(-7)$.
\end{examples}

By combining Theorem \ref{almostcompleteintcase} with Theorem
\ref{genericsplittypeforsyz} we see that almost complete
intersections in $K[X,Y,Z]$ are now well understood with respect
to the Weak Lefschetz property.

\bibliographystyle{amsalpha}


\end{document}